\documentclass[11pt]{article} 
 
\usepackage{latexsym} 
\usepackage{amsfonts} 
\usepackage{amsmath} 
\usepackage{amsthm} 
\usepackage{mathrsfs} 
\usepackage{dsfont} 
\usepackage{bbold} 
\usepackage[english]{babel} 
\usepackage{caption} 
\usepackage{epsfig} 
\usepackage{float} 
\usepackage{color} 
\usepackage{subfigure} 
\usepackage{psfrag} 
\usepackage{graphicx} 
\usepackage{epsfig} 
\usepackage{amsbsy} 
\usepackage{hyperref}

\newtheorem{theorem}{Theorem} 
\newtheorem*{prop*}{Theorem} 
\newtheorem{thm}[theorem]{Theorem} 
\newtheorem{coro}[theorem]{Corollary} 
 
\newtheorem{lemma}[theorem]{Lemma} 
 
\newtheorem{rmk}[theorem]{Remark}

\newcommand{\zerarcounters}{\setcounter{equation}{0}\setcounter{theorem}{0}} 
 
\newcommand{\beq}{\begin{equation}}
\newcommand{\eeq}{\end{equation}}

\newcommand{\ZZZ}{\mathds{Z}}

\newcommand{\RRR}{\mathds{R}} 
\newcommand{\TTT}{\mathds{T}}

\newcommand{\calK}{{\mathcal K}} 
 
\newcommand{\MM}{{\mathcal M}}

\newcommand{\calP}{{\mathcal P}}

\newcommand{\ol}{\overline} 
 
\newcommand{\Fullbox}{{\rule{2.0mm}{2.0mm}}} 
 
\newcommand{\EP}{\hfill\Fullbox\vspace{0.2cm}} 
\newcommand{\prova}{\noindent{\it Proof. }} 
\newcommand{\io}{\infty} 
\newcommand{\e}{\varepsilon} 
\newcommand{\eps}{\varepsilon} 
\newcommand{\al}{\alpha} 
\newcommand{\de}{\delta} 
\newcommand{\be}{\beta} 
\newcommand{\n}{\nu} 
 
\newcommand{\x}{\xi}

\newcommand{\f}{\varphi} 
\newcommand{\s}{\sigma} 
 
\newcommand{\del}{\partial}

\newcommand{\avg}[1]{\langle #1 \rangle}

\newcommand{\PP}{\mathtt{P}}

\newcommand{\ii}{{\rm i}}

\def\tilde#1{\widetilde{#1}}
 
\def\ins#1#2#3{\vbox to0pt{\kern-#2 \hbox{\kern#1 #3}\vss}\nointerlineskip}

\begin{document}
 
\title{\bf Locally integrable non-Liouville analytic geodesic flows}

\author{Livia Corsi$^{1}$ and Vadim Kaloshin$^{2}$
\vspace{2mm} 
\\ \small
$^{1}$ School of Mathematics, Georgia Institute of Technology, 686 Cherry St. NW, Atlanta GA, 30332, USA
\\ \small 
$^{2}$ Department of Mathematics, ETH - Inst. f\"ur Theoretische Studien,
Clausiusstrasse 47, 8092 Z\"urich, CH
\\ \small 
E-mail:  lcorsi6@math.gatech.edu, vadim.kaloshin@gmail.com}
\date{} 
 
\maketitle

\tableofcontents  


\zerarcounters
\section{Introduction and statement of the results}

Let $\MM$ be a $n$-dimensional smooth manifold with Riemanninan 
metric $\tilde{g}=(\tilde{g}_{ij})$. Using the standard transformation $p_j=\tilde{g}_{ij}\dot{q}_i$
the geodesic flow associated with $\tilde{g}$ can be regarded as an hamiltonian system

\begin{equation}\label{flusso}
\left\{
\begin{aligned}
&\dot{p} = \del_q \tilde{H} \\
&\dot{q} = - \del_p \tilde{H}
\end{aligned}
\right.
\end{equation}
where
$$
\tilde{H}(q,p) = \frac{1}{2}\sum_{i,j=1}^n \tilde{g}^{ij}(q)p_ip_j = \frac{1}{2}\sum_{i,j=1}^n \tilde{g}_{ij}(q)\dot{q}_i\dot{q}_j
$$

The geodesic flow \eqref{flusso} is called \emph{completely Liouville integrable} if it admits $n$ smooth independent functions
$H_1(q,p),\ldots,H_n(q,p)$ such that

\begin{itemize}

\item each $H_i$ is an integral of the geodesic flow, i.e. it is constant along each geodesic line $(q(t),p(t))$.

\item the functions $H_i$ Poisson-commute on $T^*\MM$, i.e. $\{H_i,H_j\}:=\del_p H_i\cdot \del_q H_j - \del_qH_i\cdot\del_pH_j=0$

\end{itemize}

Let us focus on the case $\MM=\TTT^2$ and use coordinates $q=(q_1,q_2)\in \TTT$.
If the Hamiltonian $\tilde{H}$ is of the form
$$
\tilde{H}(q,p) = (g_1(q_1) + g_2(q_2))(p_1^2+p_2^2)
$$
then the corresponding metric is said to be
 \emph{separable} or \emph{Liouville}: it is well known (see for instance \cite{Gug}) that
{every surface of revolution admits a Liouville metric}.

Suppose that the metric is ``diagonal'', namely the corresponding Hamiltonian is of the form
$$
\tilde{H}(q,p) =V(q)(p_1^2+p_2^2).
$$

We can rewrite $\tilde{g}$ as a Jacobi metric ${g}:=(e-f(q))\tilde{g}$ with
$$
e:=\max_{q\in\TTT^2} V(q),\qquad f(q) := e - V(q)
$$
so that the geodesic flow corresponds to the hamiltonian flow associated with the mechanical Hamiltonian
\begin{equation}\nonumber
H(q,p)=\frac{1}{2}|p|^2 + f(q)
\end{equation}
and the metric is separable if the Hamiltonian above is separable, namely if
the potential $f$ can be written as the sum of a function of $q_1$ only plus a function of $q_2$ only;
we refer the reader to \cite{Sch} for a recent review on the topic and some related open questions.

In \cite{BMF} the Authors prove that if a metric on $\TTT^2$ is such that the geodesic flow admits an
integral which is quadratic in the momenta then the metric is Liouville;
in  \cite{La} the Author says that Liouville metrics are the largest known class of integrable metrics.

\vskip 0.1in 
 
 \noindent
 {\bf A Forklore Conjecture } {\it If a metric on $\TTT^2$ is integrable, then it is Liouville. }

\vskip 0.1in 

The present paper provides a counterexample to the conjecture above; let us now state our result precisely.

Denote by $\calP$ the cone in the action space
\begin{equation}\label{pcone}
\calP :=\{p\in \RRR^2\;:\;-\sqrt{2}p_2<p_1<\sqrt{2}p_2,\ p_1>0\}\,.
\end{equation}
The choice of the aperture of the cone is made so that the boundaries are Diophantine directions.
It will be usefull to use polar coordinates to describe the $p$-space (see Remark \ref{fico}) i.e. we
may write
\begin{equation}\label{polare}
\calP =\{p=(\f_p,r_p)\;:\;-\arctan(\frac{1}{\sqrt{2}})<\f_p<\arctan(\frac{1}{\sqrt{2}}),\ r_p>0\}\,.
\end{equation}

Denote by $|\cdot|$ the euclidean norm of a two dimensional vector.
Our main result is the following.
\begin{thm}\label{main}
 There exists a real-on-real analytic mechanical
Hamiltonian
\begin{equation}\label{ham}
H_\e(q,p) = \frac{|p|^2}{2}+ f(q;\e) = \frac{|p|^2}{2}+\sum_{s\ge1}\e^s f_s(q)\,,
\end{equation} 
with a nonzero potential $f(q,\e)$ and an analytic change of variables $\Phi$ 
such that $H_\e\circ \Phi=|p|^2/2$ on the energy surface $\{H_\e=1/2\}$ and 
$p \in \calP$.
\end{thm}

\begin{coro}
There is a non-Liouville analytic metric on $\TTT^2$ which is 
integrable in an open set of the energy surface $\{H_\e=1/2\}$. 
 \end{coro}

Of course if one wants an example on $\TTT^n$, one can for instance decompose
$\TTT^n=\TTT^2\times\TTT^{n-2}$ and consider a metric which is the product of the metric provided by Theorem
\ref{main} for $\TTT^2$ and any integrable metric for $\TTT^{n-2}$.
However our construction strongly depend on the dimension; see Section \ref{finale} for additional comments.


On the other hand, if we allow $f$ to depend also on $p$, it is much easier to find non-separable Hamiltonian which is integrable
on a whole domain of the phase space, and moreover the construction holds in any dimension. 
Actually, denoting by $\PP$ any domain in the action space, one has the following result.

\begin{thm}\label{thm-p}
There exists a real-on-real  analytic 
Hamiltonian
\begin{equation}\label{ham-p}
H_\e(q,p) = \frac{|p|^2}{2}+ f(q,p;\e) = \frac{|p|^2}{2}+\sum_{s\ge1}\e^s f_s(q,p)\,,
\end{equation} 
with a nonzero  $f(q,p,\e)$ and an analytic change of variables $\Phi$ 
such that $H_\e\circ \Phi=|p|^2/2$ on a domain $\TTT^n\times\PP$.
\end{thm}

The proof of Theorem \ref{main} is an explicit iterative construction of the potential $f$ as convergent
power series. In Section \ref{formale} we perform by hand the first $3$ steps in order to understand the
general picture: as it is quite common in perturbation theory, from the $4$th step on the 
contruction ``stabilizes'' and in Section \ref{secs} we describe the generic step $s$ and show
the convergence of the series.
Then in Section \ref{azione} we prove Theorem \ref{thm-p}, which is rather straightforward.
 Finally in Section \ref{finale} we make some further comment about the results and their proofs.

\smallskip

\noindent
{\bf Acknowledgements}. Part of this research was performed during a period when L.C. was supported by a
CRC Postdoctoral Fellowship at McMaster University. V.K. acknowledges a partial support of the NSF grant DMS-1402164.
L.C. acknowledges a partial support of the NSF grant DMS-1500943.

\zerarcounters
\section{An iterative procedure}\label{formale}

 Our aim is to explicitely construct a real-on-real potential 
 \begin{equation}\label{poto}
 f(q;\e) = \sum_{s\ge1}\e^s f_s(q)
 \end{equation}
 and a change of variables $\Phi$ such that $H_\e\circ \Phi = |p|^2/2 + h(p;\e)$ on the level surface $\{H_\e = 1/2\}$,
 for some
 \begin{equation}\label{acca}
 h(p;\e) = \sum_{s\ge1} \e^s h_s(p)\,.
\end{equation}
 
As usual in KAM-like problems, the change of variables $\Phi$ will be the time-1 map
generated by an Hamiltonian
\begin{equation}\label{G}
G(p,q)=G(p,q;\e)=\sum_{s\ge1}\e^s G_s(p,q)\,,
\end{equation}
so that
\begin{equation}\label{poisson}
H\circ \Phi = \sum_{n\ge0}\frac{1}{n!}\{H,G\}^{(n)}\,,
\end{equation}
where we used the notation
\begin{equation}\label{nota}
\begin{aligned}
&\{H,G\}^{(0)}:= H\,,\\
&\{H,G\}^{(n)}:=\{\{H,G\}^{(n-1)},G\}\,;
\end{aligned}
\end{equation}
see for instance \cite{CKO}.
 
We shall construct $f$ and $\Phi$ via an iterative procedure: at each step we suitably fix $f_s$ and $G_s$, and
we provide appropriate bounds from which we eventually infer the convergence of both the series \eqref{poto} and \eqref{G}.

Let $\calK$ be the dual cone
\begin{equation}\label{kcone}
\calK:=\{k\in \ZZZ^2\setminus\{0\}\;: |k\cdot p| \ge \frac{1}{2}|p||k|,\mbox{ for all }p\in\calP\}\,.\ 
\end{equation}
 
We start our procedure by chosing any non-separable $f_1(q)$ whose Fourier modes are supported in $\calK$, i.e.
\begin{equation}\label{f1}
f_{1}(q)=\sum_{k\in \calK} f_{1,k}e^{\ii k\cdot q}\,.
\end{equation}

Of course, by analyticity we have
\begin{equation}\label{anal}
|f_{1,k}|\le M e^{-\x_0 |k|}
\end{equation}
for some given positive constants $M,\x_0$. Moreover, in order for $f_1$ to be real-on-real, we need to require
\begin{equation}\label{realonreal1}
f_{1,-k} = \ol{f}_{1,k},
\end{equation}
where for a complex number $z$ we denoted by $\ol{z}$ its complex conjugate.


Collecting together the same orders in $\e$ and denoting 
$H_0=H_0(p ):=|p|^2/2$ we have
\begin{equation}\label{expand}
H\circ \Phi = H_0(p )+\sum_{s\ge1}\e^{s}\Big(f_s+
\sum_{m\ge1}^s
\frac{1}{m!}\sum_{\substack{n_0+\ldots+n_m=s\\ n_0\ge0\\n_1,\ldots,n_m\ge1}}\{\{\ldots\{\{f_{n_0},G_{n_1}\},
G_{n_2}\}\ldots\},G_{n_m}\}
\Big)\,.
\end{equation}

For all ${\mathtt n}\ge1$, we denote
\begin{equation}\label{tronco}
H^{({\mathtt n})}_\e=H_0(p )+\sum_{s=1}^{\mathtt n} \e^{s}f_s(q)\,,
\end{equation}
with $f_1(q)$ given in \eqref{f1} 
while $f_s(q)$ for $s\ge 2$ are still to be found.

We now compute explicitly the first three orders in order to understand the general behavior.


 \subsection{The first and second steps}

At order $\e$, the r.h.s. of \eqref{expand} reads
\begin{equation}\label{uno}
f_1+\{H_0,G_1\} = f_1 + p\cdot\del_q G_1\,,
\end{equation}
so that, defining formally
\begin{equation}\label{G1}
G_1(p,q) = - \sum_{k\in \calK}\frac{ f_{1,k}}{\ii p\cdot k} e^{\ii k\cdot q}
\end{equation}
we have
\begin{equation}\label{primo}
\{H_0,G_1\}=-f_1\,.
\end{equation}
Note that, due to \eqref{realonreal1}, $G_1$ is real-on-real as well.

For $1/2<|p|<2$ we have the bound
\begin{equation}\label{int}
|G_{1,k}|\le \frac{2M}{|k|}e^{-\x_0 |k|}\,,
\end{equation}
so that the function $G_1(p,q)$ in \eqref{G1} is well defined on $(\{1/2<|p|<2\}\cap\calP)\times \TTT^2$.
Moreover we can take $h_1(p)\equiv 0$.

At order $\e^2$ the r.h.s. of \eqref{expand} is
\begin{equation}\label{due}
f_2+\{f_1,G_1\}+\{H_0,G_2\}+\frac{1}{2}\{\{f_0,G_1\},G_1\} \stackrel{\eqref{primo}}{=} 
f_2+\frac{1}{2}\{f_1,G_1\}+\{H_0,G_2\}\,.
\end{equation}
We start by imposing {that on the energy surface
$\{H_0=1/2\}$ we have }
\begin{equation}\label{secondo}
\{H_0,G_2\}=-(f_2+\frac{1}{2}\{f_1,G_1\}) + h_2\,,
\end{equation}
which in Fourier reads
\begin{subequations}\label{homof}
\begin{align}
&(\ii p\cdot k)G_{2,k}=-(f_{2,k}+(\frac{1}{2}\{f_1,G_1\})_k)\,,\qquad k\ne 0
\label{range} \\
& h_2 =f_{2,0}+(\frac{1}{2}\{f_1,G_1\})_0\,,\qquad k= 0.
\label{compa}
\end{align}
\end{subequations}
Note that \eqref{compa} leaves us $f_{2,0}$ as a free parameter while
 we need to define $f_{2,k}$ is such a way that the r.h.s. of 
\eqref{range} is zero when $|p\cdot k|$ is {$O(\e)$} small.
{We can't do it in a uniform way, so we solve it with 
a precision $O(\e)$}.
 
For $k\ne 0$, denote
\begin{equation}\label{perp}
p^0_k:= k^{\perp}/|k|\,,
\end{equation}
and note that the Poisson bracket appearing in \eqref{homof} is non-zero only for $k\in\calK_2$ defined as
\begin{equation}\label{doppiocono}
\calK_2:= \{k\in\ZZZ^2\;:\; k = k_1 + k_2\,, \mbox{ for some }k_1,k_2\in\calK\}.
\end{equation}

We distinguish two subset of $\calK_2$, namely
\begin{equation}\label{k2grande}
\calK_2^{big} := \{k\in\calK_2\;:\; p^0_k\notin\calP\}
\end{equation}
and 
\begin{equation}\label{k2piccolo}
\calK_2^{small} := \{k\in\calK_2\;:\; p^0_k\in\calP\},
\end{equation}
and we analyze the Fourier modes separately.

\noindent
{\bf case 1.} $k\in \calK_2^{big}$. In this case the l.h.s. of \eqref{homof} cannot vanish.
Of course, although $p^0_k\notin\calP$ there might be $p\in\calP$ such that $p\cdot k$ is ``too small'': this might happen if $p^0_k$ is
``close'' to the boundary of $\calP$, so we set
\begin{equation}\label{f2bordo}
f_{2,k}:=-\frac{1}{2}\sum_{k_1+k_2=k}
\frac{k_1 f_{1,k_1} k_2 f_{1,k_2}}{(\ol{p}_k\cdot k_2)^2}\,,
\end{equation}
where $\ol{p}_k$ is the minimizer of $|p - p^0_k|$ for $p$ varying in the closure of $\calP\cap\{|p|=1\}$, i.e. it is either
 $(\frac{\sqrt{2}}{3},\frac{1}{3})$ or $(\frac{\sqrt{2}}{3},-\frac{1}{3})$.
Note that in both cases $\ol{p}_k$ is a Diophantine vector so we do not have to worry about the smallness of denominators of the form $\ol{p}_k\cdot k'$
for any $k'\in\ZZZ^2\setminus\{0\}$.

\medskip

\noindent
{\bf case 2.}  $k\in\calK_2^{small}$. In this case the l.h.s. of \eqref{homof}  vanishes when $p=p^0_k$ so that we first need to impose
\begin{equation}\label{f2}
f_{2,k}:=-\frac{1}{2}\sum_{k_1+k_2=k}
\frac{k_1 f_{1,k_1} k_2 f_{1,k_2}}{(p^0_k\cdot k_2)^2}\,.
\end{equation}

\medskip

In both cases 1,2 we can define, at least formally,
\begin{equation}\nonumber
G_{2,k}=\frac{-1}{\ii p\cdot k}
\Big(f_{2,k}+\frac{1}{2}\sum_{k_1+k_2=k}
\frac{-k_1 f_{1,k_1} k_2 f_{1,k_2}}{(p\cdot k_2)^2}
\Big)\,,
\end{equation}
and 
\begin{equation}\label{G2}
{G^0_{2,k}}=\frac{-1}{\ii p\cdot k}
\Big(f_{2,k}+ \frac{{{|p|^2}}}{2}\sum_{k_1+k_2=k}\,
\frac{-k_1 f_{1,k_1} k_2 f_{1,k_2}}{(p\cdot k_2)^2}
\Big)\,.
\end{equation}

\begin{rmk}\label{diverso}
When $|p|^2-1=O(\e)$, $G_2$ and $G_2^0$ differ by $O(\e)$
and the mismatch goes to the next order. In particular, $G_2$ 
will turn out to be not well defined in this neighbourhood, but $p\cdot\partial_q G_2$ is well defined there.
\end{rmk}

\begin{rmk}\label{allafaccia}
Note that in the sum apperaring in \eqref{f2} in principle the term with $k_1=k_2$ might be the source of a problem.
Indeed if $k$ is such that there exists $k_1$ so that $2k_1 = k$, then the denominator in \eqref{f2} satisfies
$$
p^0_k\cdot k_2 = p^0_k\cdot k_1 = \frac{1}{2} p^0_k\cdot k =0.
$$
However this cannot happen because in this case $k\parallel k_1$, which implies $k\in\calK\subseteq\calK_2^{big}$
so that $f_{2,k}$ is given by \eqref{f2bordo} and not by \eqref{f2}.
\end{rmk}

\begin{rmk}\label{realonreal2}
Note that, because of \eqref{realonreal1}, in both cases 1 and 2 one has
\begin{equation}\nonumber
\begin{aligned}
f_{2,-k} &= -\frac{1}{2}\sum_{k_1+k_2=-k}
\frac{k_1 f_{1,k_1} k_2 f_{1,k_2}}{(\tilde{p}\cdot k_2)^2} \\
&= -\frac{1}{2}\sum_{k_1+k_2=k}
\frac{(-k_1) f_{1,-k_1} (-k_2) f_{1,-k_2}}{(-\tilde{p}\cdot k_2)^2} \\
&=-\frac{1}{2}\sum_{k_1+k_2=k}
\frac{k_1 \ol{ f}_{1,k_1} k_2 \ol{f}_{1,k_2}}{(\tilde{p}\cdot k_2)^2} = \ol{f}_{2,k}
\end{aligned}
\end{equation}
where $\tilde{p}$ is equal to $\ol{p}_k$ in case 1 and to $p^0_k$ in case 2.
\end{rmk}

Regarding the bounds, first of all we see that in both cases 1 and 2 one has
\begin{equation}\label{stimof2}
|f_{2,k}|\le \frac{M^2}{2}\sum_{k_1+k_2=k}\frac{|k_1|}{|k_2|}e^{-\x_0(|k_1|+|k_2|)} 
\le M^2 e^{-\x_0 |k|/2}\,,
\end{equation}
then we notice that, setting
\begin{equation}\label{F2}
F_{2,k}(p ):=\frac{{{|p|^2}}}{2}\sum_{k_1+k_2=k}\frac{-k_1 f_{1,k_1} k_2 f_{1,k_2}}{(p\cdot k_2)^2}\,,
\end{equation}
 in case 1 one has
\begin{equation}\label{increfacile}
f_{2,k}+\frac{1}{2}\{f_1,G_1\}_k =  F_{2,k} (\ol{p}) - F_{2,k}(p )
\end{equation}
whereas in case 2 we have
\begin{equation}\label{incre}
f_{2,k}+\frac{1}{2}\{f_1,G_1\}_k =  F_{2,k} (p_k^0) - F_{2,k}(p ).
\end{equation}

Note that clearly, as in Remark \ref{realonreal2} one has
\begin{equation}\label{roar}
F_{2,-k} = \ol{F}_{2,k}.
\end{equation}

\begin{rmk}\label{guarda}
One has
$$
p\cdot\del_p F_{2,k}(p ) = {{{|p|^2}}}\sum_{k_1+k_2=k}\frac{-k_1 f_{1,k_1} k_2 f_{1,k_2}}{(p\cdot k_2)^2} +
|p|^2\sum_{k_1+k_2=k}\frac{k_1 f_{1,k_1} k_2 f_{1,k_2}}{(p\cdot k_2)^2} =0\,.
$$
\end{rmk}

Given two vectors $u,v\in\RRR^2$ we denote by $\f(u,v)$ the smaller angle between the two vectors,
so with this notation we can rewrite
\begin{equation}\label{F2meglio}
F_{2,k}(p ):=\frac{1}{2}\sum_{k_1+k_2=k}\frac{-k_1 f_{1,k_1} k_2 f_{1,k_2}}{| k_2|^2\cos^2(\f(p,k_2))}\,.
\end{equation}

\begin{rmk}\label{fico}
Note that if $p,p'$ are parallel, then 
$F_{2,k}(p )-F_{2,k}(p' ) =0$. This is the reason why it is convenient to describe the $p$-variables
in polar coordinates as in \eqref{polare}: indeed with that notation we have
$F_{2,k}(p)=\ol{F}_{2,k}(\f_p)$, i.e. it is a function of the angular variable only.
\end{rmk}

Using the notation \eqref{F2meglio} and setting
\begin{equation}\label{ptilde}
\tilde{p}_k := \left\{
\begin{aligned}
&\ol{p}_k\qquad k\in \calK_2^{big} \\
&p^0_k\qquad k\in \calK_2^{small}
\end{aligned}
\right.
\qquad\qquad
\hat{p}_k :=
 \left\{
\begin{aligned}
&0\qquad k\in \calK_2^{big} \\
&p^0_k\qquad k\in \calK_2^{small}
\end{aligned}
\right.
\end{equation}
we see that 
\begin{equation}\label{eccola}
G^0_{2,k} = \frac{F_{2,k}(p ) - F_{2,k} (\tilde{p}_k)}{\ii (p- \hat{p}_k)\cdot k}
\end{equation}
and hence
$$
G^0_{2,-k} = \ol{ G^{0}}_{2,k}.
$$
Moreover we can bound
\begin{equation}\label{deriva}
\begin{aligned}
|{G^0_{2,k}}| &= \frac{|F_{2,k}(p ) - F_{2,k} (\tilde{p}_k)|}{|(p- \tilde{p}_k)\cdot k|}\\
&\le
\frac{1}{2}\sum_{k_1+k_2=k}
\frac{|k_1 f_{1,k_1} k_2 f_{1,k_2}|}{| k_2|^2\cos^2(\f(p,k_2))\cos^2(\f(\tilde{p}_k,k_2))}
\frac{|\cos^2(\f(p,k_2))-\cos^2(\f(\tilde{p}_k,k_2))|}{|(p-\tilde{p}_k)\cdot k|}\,.
\end{aligned}
\end{equation}

Now if $|(p-\tilde{p}_k)\cdot k| =\de$, then $\cos(\f(p,k))=(|k|r_p)^{-1}\de$ which in turn implies
$$
|\f(p,k)-\f(\tilde{p}_k,k)|\le \frac{2\de}{|k|r_p}\,,
$$
and hence
\begin{equation}\label{cos}
|\cos^2(\f(p,k_2))-\cos^2(\f(\tilde{p}_k,k_2))| \le 4\left|\sin\big(\frac{\f(p,k_2)-\f(\tilde{p}_k,k_2)}{2}\big)\right|
\le \frac{2\de}{|k|r_p}\,,
\end{equation}
from which we deduce
\begin{equation}\label{finita}
\begin{aligned}
|{G^0_{2,k}}| 
&\le \sum_{k_1+k_2=k}
\frac{2|k_1 f_{1,k_1} k_2 f_{1,k_2}|}{|k|^2r_p^2| k_2|^2\cos^2(\f(p,k_2))\cos^2(\f(\tilde{p}_k,k_2))}
\le CM^2e^{-3\x_0|k|/4}\,,
\end{aligned}
\end{equation}
for some positive constant $C$,
where in the last inequality we used \eqref{stimof2} and the fact that $r_p^2=|p|^2=1+O({\e})\le 2$.
Note that we are left with the average $G^0_{2,0}$ as a further free parameter.
Note also that $G^0_2$ is analytic as function of $p$ in a $O(\sqrt \e)$-neighbourhood of 
$\{H_0=1/2\}$.

In conclusion $G^0_2$ is a real-on-real function, well defined and analytic in a $O(\sqrt \e)$-neighbourhood of 
$\{H_0=1/2\}$ and we have the uniform bound \eqref{deriva} for its Fourier coefficients;
on the other hand $G_2$ is well defined and its Fourier 
coefficients admit the bound \eqref{deriva} {only on the surface $\{H_0=1/2\}$},
being equal to $G^0_2$ on such surface. Moreover the averages $f_{2,0}$ and $G^0_{2,0}$ are still free parameters.


{
\subsection{The $3$-rd step}

At the previous step we obtained
\begin{equation}\label{aol}
f_2(q) = \sum_{k\in \ZZZ^d} f_{2,k}e^{\ii k\cdot q}
\end{equation}
with $f_{2,0}$ a free parameter and $f_{2,k}$
 given by \eqref{f2bordo} or (\ref{f2}) according on $k$, as the restriction of 
$\{f_1,G_1\}$ to the  surface $\{H_0=1/2\}$.

We now refine this definition in order to fit the restriction onto the  energy 
surface $\{H_\e^{(1)}=1/2\}$. In other words we want the r.h.s. of (\ref{secondo}) to
vanish on $\{H_\e^{(1)}=1/2\}$. Let $p^1_k(q,\e)$ be either $\ol{p}_k$ for $k\in\calK_2^{big}$ or the unique point in $\calP$ such that
\begin{equation}\label{inter}
\begin{aligned}
p\cdot k =0\qquad \qquad \\
\frac{1}{2}|p|^2 + \e f_1(q) = \frac{1}{2}\,,
\end{aligned}
\end{equation}
for $k\in\calK_2^{small}$.
Notice that for $k\in\calK_2^{small}$, $p_k^1(q,\e)$ can be written in the form 
$$
p_k^1(q,\e)=p_k^0(1-2 \e f_1(q))^{1/2}\,;
$$ 
Now set
\begin{equation}\label{f21}
{f}^1_{2,k}(q,\e):=-\frac{1}{2}\sum_{k_1+k_2=k}\frac{k_1 f_{1,k_1} k_2 f_{1,k_2}}{(p^1_k(q,\eps)\cdot k_2)^2}\,,
\end{equation}
%
%
and note that
$$
{f}^1_{2, -k}(q,\e) = \ol{{f}^1}_{2,k}(q,\e)\,.
$$

Define also
\begin{equation}\label{3star}
\begin{aligned}
f_3^*(q)&:=\lim_{\e \to 0}\ \sum_k  \left(
\frac{{f}^1_{2,k}(q,\e)-f_{2,k}}{\eps} -
\frac{\left.\{H_0,G_{2}^0\}_k\right|_{p=p_k^1(q,\e)}}{\eps}\right) e^{ \ii k\cdot q}\\
&= \lim_{\e \to 0}\  \sum_k  \left(
\frac{{f}^1_{2,k}(q,\e)-f_{2,k}}{\e}-
\frac{i p\cdot kG_{2,k}^0(p_k^1(q,\e))}{\eps} \right) e^{ \ii k\cdot q}\\
&= f_1(q)\Big(f_2(q)+\sum_k p_k^0\cdot \del_p F_{2,k}(p_k^0)e^{\ii k\cdot q}\Big) \stackrel{{\rm Rmk}\,\ref{guarda}}{=}
f_1(q)f_2(q)\,,
\end{aligned}
\end{equation}
where we used also the fact that
$$
\lim_{\e\to0}p\cdot\partial_qG_2^0(q,p_k^1(q,\e))=
p\cdot\del_q G_{2,k}(q,p_k^0)=0
$$
by construction. In particular we see that $f^*_3(q)$ is also real-on-real.

 Using the formal definitions of $G_1,G_2$, the equation for $G_3$ formally reads
\begin{equation}\label{G3fo}
\begin{aligned}
p\cdot \del_q G_3 &=h_3 -\Big(
f_3+{f_3^*}
+\{f_1,G_2\}+\{f_2,G_1\} \\
&\qquad+\frac{1}{2}\left(\{ \{H_0,G_2\},G_1\}
+ \{ \{H_0,G_1\},G_2\}\right) +\frac 16 \{ \{\{H_0,G_1\},G_1\},G_1\}
\Big) \\
&=h_3{-}\Big(
f_3+{f_3^*}
+\frac{1}{2}\{f_1,G_2\}+\frac{1}{2}\{f_2,G_1\}  +\frac{1}{2} \{ \{\{f_1,G_1\},G_1\}
\Big)
\end{aligned}
\end{equation}
namely in Fourier, again formally,
\begin{subequations}\label{G3ff}
\begin{align}
&(\ii p \cdot k)G_{3,k}={-}\Big( f_{3,k}+{f^*_{3,k}}+\frac{1}{2}\{f_1,G_2\}_k+\frac{1}{2}\{f_2,G_1\}_k
+\frac{1}{12}  \{ \{f_1,G_1\},G_1\}_k
\Big),\quad k\ne0 
\label{range} \\
&h_3=\Big( f_{3,0}+{f^*_{3,0}}+\frac{1}{2}\{f_1,G_2\}_0+\frac{1}{2}\{f_2,G_1\}_0
+\frac{1}{12}  \{ \{f_1,G_1\},G_1\}_0
\end{align}
\end{subequations}
where $\{\{\cdots\},\cdot\}_k$ is the $k$-th
Fourier coefficient of the underlying function.} 
Notice that we have ``energy reduction'' correction term
coming from the difference between restricting $\{f_1,G_1\}$ 
to $\{H_0=1/2\}$ and to $\{H^{(1)}_\e=1/2\}$.

As before we note that the Poisson brackets above are non-zero only for
\begin{equation}\label{k3}
k\in\calK_3 := \{k\in\ZZZ^2\;:\; k = k_1 + k_2 + k_3 \mbox{ for some }k_1,k_2,k_3\in\calK\},
\end{equation}
and again we distinguish the two subsets
\begin{equation}\label{k3grande}
\calK_3^{big} := \{k\in\calK_3\;:\; p^0_k\notin\calP\}
\end{equation}
and 
\begin{equation}\label{k3piccolo}
\calK_3^{small} := \{k\in\calK_3\;:\; p^0_k\in\calP\},
\end{equation}
and as before we define the $k\ne0$ Fourier coefficients $f_{3,k}$ differently for modes in $\calK_3^{big}$ or in $\calK_3^{small}$

\noindent
{\bf case 1.} $k\in \calK_3^{big}$.  In this case we set (recall that $\ol{p}_k$ is Diophantine)
\begin{equation}\label{f3bordo}
f_{3,k}:= -f^*_{3,k}
-\frac{1}{2}\left.\{f_1,G_2^0 \}_k\right|_{p=\ol{p}_k}
-\frac{1}{2}\left.\{f_2,G_1\}_k\right|_{ p=\ol{p}_k}
-\frac{1}{12}
\left.\{\{f_1,G_1 \},G_1 \}_k\right|_{p=\ol{p}_k}
\end{equation}

\medskip

\noindent
{\bf case 2.} $k\in \calK_3^{small}$. In this case we set
\begin{equation}\label{f3}
f_{3,k}:= -f^*_{3,k}
-\frac{1}{2}\left.\{f_1,G_2^0 \}_k\right|_{p={p}^0_k}
-\frac{1}{2}\left.\{f_2,G_1\}_k\right|_{ p={p}^0_k}
-\frac{1}{12}
\left.\{\{f_1,G_1 \},G_1 \}_k\right|_{p={p}^0_k}.
\end{equation}

Note that one can reason as in Remark \ref{allafaccia} to deduce that in \eqref{f3} above no zero divisor appear.
Notice also that by construction $f_{3}$ is real-on-real.

{Clearly \eqref{G3fo} and \eqref{G3ff} are just formal expressions since $G_2$ is defined only on the surface
$\{H_0=1/2\}$. So we proceed as in the previous case,}
i.e. we modify the definition of $G_{3,k}$ obtained by solving \eqref{G3fo} (once $f_{3,k}$ is defined according either to
\eqref{f3bordo} or \eqref{f3}) to 
\begin{equation}\label{G30}
\begin{aligned}
{G_{3,k}^0}(p)=
\frac{-1}{\ii p\cdot k}\Big( f_{3,k}+{f^*_{3,k}}+
\frac{|p|^2}{2}\{f_1,G_2^0\}_k+\frac{|p|^2}{2}\{f_2,G_1\}_k
+\frac{|p|^4}{12}  \{ \{f_1,G_1\},G_1\}_k
\Big)\,.
\end{aligned}
\end{equation}
Again by definition one has
$$
G_{3,-k}^0 (p) = \ol{G^0}_{3,k} (p).
$$

{Notice that in \eqref{G30} we use $G_2^0$, which is well defined and uniformly 
bounded in a $O(\sqrt{\e})$-neighbourhood of $\{H_0=1/2\}$.} 
As done for $G_2^0$ we obtain a uniform bound for $G^0_3$ in a $O(\e)$-neighbourhood of $\{H_\e^{(1)}=1/2\}$.
Indeed setting
\begin{equation}\label{F3}
\begin{aligned}
F_{3,k}(p ):=  \frac{|p|^2}{2}\{f_1,G_2^0\}_k+\frac{|p|^2}{2}\{f_2,G_1\}_k
+\frac{|p|^4}{12}  \{ \{f_1,G_1\},G_1\}_k\,,
\end{aligned}
\end{equation}
we see that
$$
|G^0_{3,k}(p) | = \frac{|F_{3,k}(p ) - F_{3,k} (\tilde{p}_k)|}{|(p- \hat{p}_k)\cdot k|}\,,
$$
where we are using the notation \eqref{ptilde} with $\calK_2\rightsquigarrow\calK_3$.
 Hence we can reason exactly as in \eqref{deriva} to get $|G^0_{3,k}(p )|\le CM^3e^{-5\x_0|k|/8}$.

\begin{rmk}\label{normo}
As for the second step, we introduced the factors $|p|^a$, $a=2,4$, in order to make sure that $F_{3,k}$ depend
on $p$ only through $\f_p$, and this is the reason behind the choice of the exponents; see also Remark
\ref{fico}.
Indeed from \eqref{eccola} we see that $G^0_{2,k}$ is of the form
$$
G^0_{2,k} = \frac{1}{p\cdot k}g_k(\f_p)
$$
for some function $g_k$ depending on $p$ only through $\f_p$, and thus by explicit computation one sees that $|p|^2\{f_1,G_2^0\}_k$
depends on $p$ only through $\f_p$; the same type of argument apply for the other  terms.
In other words, the normalization factors $|p|^a$ in \eqref{F3} are made so that 
$F_{3,k}(p)=\ol{F}_{3,k}(\f_p)$.
\end{rmk}

Let {$G^{(3)}:=\e G_1+\eps^2 G^0_2+\eps^3 G^0_3$} which is real-on-real by construction, and $\Phi^{(3)}$
be the time-$1$ map of $G^{(3)}$. By construction, fixing $f_s(q)$ as above for $s=1,2,3$
while $f_s(q)$ are still arbitrary for $s\ge4$, we get
\[
H \circ \Phi^{(3)}= \frac{p^2}{2}+\e h_1(p) +\e^2 h_2(p) + \e^3 h_3(p) + \eps^4 R_4\,,
\]
where $R_4$ is a suitable remainder; precisely
\begin{equation}\label{R4}
\begin{aligned}
\e^4 R_4&= \frac{1}{2}\Big(
\{\{H,\e G_1^0\}, \e^3G^0_3\} + \{\{H,\e^2G_2^0\},\e^2G_2^0+\e^3G_3^0\}+\{\{H,\e^3 G_3^0\},G^{(3)}\}
\Big)\\
&+\frac{1}{3!}\Big(
\{\{\{H,\e G^0_1\},\e G_1^0\},\e^2G_2^0+\e^3G_3^0\} + \{\{\{H,\e G^0_1\},\e^2G_2^0+\e^3G_3^0\},G^{(3)}\}\\
&\qquad+\{\{\{H,\e^2G_2^0+\e^3G_3^0\} ,G^{(3)}\},G^{(3)}\}
\Big)\\
&+\sum_{n\ge 4}\frac{1}{n!}\{H,G^{(3)}\}^{(n)}\,.
\end{aligned}
\end{equation}
Note that each term in \eqref{R4} above has at least a factor $\e^4$.
Note also that $G^{(3)}$ is analytic as function of $p$.



{
\section{The $s$-th step}\label{secs}

We now describe the procedure at the $s$-th step: as we shall see we first need to perform
formal computations and then provide suitable modifications that allow us to obtain the bounds
needed for the convergence of the algorithm.

\subsection{Formal expansion}\label{sec.formale}

Consider the formal truncation of $H \circ \Phi$ at order $s-1$, namely
\[
H_0(p)+\e (f_1+\{H_0,G^0_1\}) +
\e^2 (f_2+\{H_0,G^0_2\}+\{f_1,G^0_1\})+\]\[
+\sum_{n=3}^{s-1}\e^n\Big(
\sum_{m=1}^{n}\frac{1}{m!}
\sum_{\substack{n_0+\ldots+n_m=n\\ n_0\ge0\\n_1,\ldots,n_m\ge1}}
\{\{\ldots\{\{f_{n_0},G^0_{n_1}\},
G^0_{n_2}\}\ldots\},G^0_{n_m}\}\Big),
\]
where we denoted $G_1^0=G_1$.
Recall that $f_2$ as in (\ref{f2}) annihilate the $O(\e^2)$
term, when restricted to the energy surface $\{H_0=1/2\}$. 
Then using (\ref{f21}) we produced a correction term $f_3^*$ 
so that $f_2+\e f_3^*$ annihilate the $O(\e^2)$
term to the order $O(\e^4)$, when restricted to the energy 
surface $\{H_\e^{(1)}=1/2\}$ and so on; eventually we achieve 
the $O(\e^{s+1})$ order cancellation on 
$\{H_\e^{(s-1)}=1/2\}$ we need to cancel
the $\e^2(\dots)$ bracket to $O(\e^{s-1})$, 
the $\e^3(\dots)$ bracket to $O(\e^{s-2})$, and so on
the $\e^{s-1}(\dots)$ bracket to $O(\e^2)$.

Assume recursively that for all $n=1,\ldots,s-1$ the functions $G^0_{n,k}$ are analytic for $p$ in a neighborhood of $\{H_\e^{(s-1)}=1/2\}\cap\calP$,
and one has
\begin{equation}\label{corsetta}
f_{n,-k} = \ol{f}_{n,k} \,, \qquad f^*_{n,-k} = \ol{f^*}_{n,k} \,, \qquad
G_{n,-k}^0 = \ol{G^0}_{n,k} \,.
\end{equation}
Of course we have to prove that \eqref{corsetta} above is satisfied also at step $s$.

For all $j=1,\ldots s$ set
\begin{equation}\label{Kj}
\calK_j=\{k\in\ZZZ^2\;:\; k=\sum_{l=1}^j k_l\; \mbox{ for some }k_1,\ldots,k_j\in\calK\}\,,
\end{equation}
distinguish the two subsets
\begin{equation}\label{kjgrande}
\calK_j^{big} := \{k\in\calK_j\;:\; p^0_k\notin\calP\}
\end{equation}
and 
\begin{equation}\label{kjpiccolo}
\calK_j^{small} := \{k\in\calK_j\;:\; p^0_k\in\calP\},
\end{equation}
and for all $j=1,\ldots s-2$ let $p^{j}_k(q,\e)$  either $\ol{p}_k$ for $k\in\calK_{s-1}^{big}$ or the unique point in $\calP$ such that
\begin{equation}\nonumber
\begin{aligned}
p\cdot k =0\qquad \qquad \\
\frac{1}{2}|p|^2 + \sum_{k=1}^{j}\e^k f_k(q) = \frac{1}{2}\,.
\end{aligned}
\end{equation}
Notice that for $k\in\calK_j^{small}$ such points can be written in the form 
\begin{equation}\label{stoqua}
p_k^{j}(q,\e)=p_k^0(1-\sum_{l=1}^{j}\e^l f_l(q))^{1/2}\,,
\end{equation}
while we recall that for $k\in\calK_j^{big}$ one has that $\ol{p}_k$ is Diophantine.

We define the  correction to 
\[
\sum_{m=1}^{n}\frac{1}{m!}
\sum_{\substack{n_0+\ldots+n_m=n\\ n_0\ge0\\n_1,\ldots,n_m\ge1}}
\{\{\ldots\{\{f_{n_0},G^0_{n_1}\},
G^0_{n_2}\}\ldots\},G^0_{n_m}\}
\]
 as
\begin{equation}\label{effepi}
\begin{aligned}
f_{s,k}^n(q)=
\lim_{\e\to0}
\frac{1}{\e^{s-n-1}} 
\sum_{m=1}^{n}\frac{1}{m!} &\Big(
\sum_{\substack{n_0+\ldots+n_m=n\\ n_0\ge0\\n_1,\ldots,n_m\ge1}}
\left.
\{\{\ldots\{\{f_{n_0}+f_{n_0}^*,G^0_{n_1}\},
G^0_{n_2}\}\ldots\},G^0_{n_m}\}\right|_{p=p_k^{s-n-1}(q,\e)}\\
&\left.
-\{\{\ldots\{\{f_{n_0} + f_{n_0}^*,G^0_{n_1}\},G^0_{n_2}\}\ldots\},G^0_{n_m}\}\right|_{p=p_k^{s-n-2}(q,\e)} \Big)\,,
\end{aligned}
\end{equation}
}
where we denoted  $f_1^*=f_2^* \equiv0$ and $f_0 = H_0$.

{
Set
\begin{equation}\label{stellina}
f_s^*(q):=\sum_{n=1}^{s-2} f_s^n(q)=
\sum_k f^*_{s,k} e^{ \ii k\cdot q}\,,
\end{equation}
and
\begin{equation}\label{margarida}
\tilde{p}_k := \left\{
\begin{aligned}
&\ol{p}_k\qquad k\in \calK_{s}^{big}\,, \\
&p^0_k\qquad k\in \calK_s^{small}\,.
\end{aligned}
\right.
\end{equation}

We define the correction term $f_s(q)$ by setting  its Fourier coefficients  as
\begin{equation}\label{effesstar}
f_{s,k}:=-\Big({f_{s,k}^*}+
\sum_{m=1}^{s}\frac{1}{m!}
\sum_{\substack{n_0+\ldots+n_m=s\\ n_0\ge0\\n_1,\ldots,n_m\ge1}}
\left.
\{\{\ldots\{\{f_{n_0}+f_{n_0}^*,G^0_{n_1}\},
G^0_{n_2}\}\ldots\},G^0_{n_m}\}\right|_{p=\tilde{p}_k}\Big)\,.
\end{equation}

Note that $f_s(q)$ is well defined since the functions $G^{0}_{n,k}$ are analytic in $p$.

\begin{rmk}\label{muscoli}
Note that by construction
$$
f_{s,-k}= \ol{f}_{s,k}\,,\qquad f^*_{s,-k}= \ol{f^*}_{s,k}
$$
so that the first two recursive assumptions in \eqref{corsetta} are satisfied also at step $s$.
\end{rmk}

\begin{rmk}\label{bada}
Both in \eqref{effepi} and \eqref{effesstar} we need to use $G^0_{n_i}$ because
the functions $G_{n_i}$ are not well defined in an open neighborhood of $|p|=1$
for $n_1\ge2$. However, their derivatives are formally
equal. 
\end{rmk}

The formal equation for $G_s$ is, therefore,
\begin{equation}\label{Gs}
p\cdot\del_q G_{s}=-\Big(
f_{s}+f^*_{s}+
\sum_{m\ge1}\frac{1}{m!}
\sum_{\substack{n_0+\ldots+n_m=s\\ n_0\ge0\\n_1,\ldots,n_m\ge1}}
\{\{\ldots\{\{f_{n_0}+f_{n_0}^*,G_{n_1}\},
G_{n_2}\}\ldots\},G_{n_m}\}\Big),
\end{equation}
but actually there are cancellations allowing to get rid of the case $n_0=0$ simply changing
the combinatorial factors, as the following result shows.

\begin{lemma}\label{lem.cancella}
The formal equation \eqref{Gs} is equivalent to
\begin{equation}\label{vera}
p\cdot \del_q G_s = -\Big(
f_s+ f_s^*+
\sum_{m=1}^{s-1}c_m
\sum_{\substack{n_0+\ldots+n_m=s\\ n_i\ge1}}\{\{\ldots\{\{f_{n_0}+f_{n_0}^*,G_{n_1}\},
G_{n_2}\}\ldots\},G_{n_m}\}
\Big)\,,
\end{equation} 
where $f_1^*=f_2^*\equiv0$ and
\begin{equation}\label{cm}
\begin{aligned}
&c_0:=1\,,\\
&c_m:=\frac{1}{m!}-\sum_{j=1}^m \frac{1}{(j+1)!}c_{m-j}\,.
\end{aligned}
\end{equation}
\end{lemma}

\prova
We prove the result by double induction on $s\ge2$ $m\ge1$. The case $s=2$, $m=1$ is the
explicit computation in \eqref{due}. Assume inductively
the statement to be true up to order $s$ with the coefficients $c_1,\ldots,c_{s-1}$ given by \eqref{cm}.
At order $s+1$ the equation \eqref{Gs} reads
\begin{equation}\label{Gs+1}
\begin{aligned}
p\cdot \del_q &G_{s+1} = -\Big(
f_{s+1}+f_{s+1}^*+
\sum_{m\ge1}\frac{1}{m!}
\sum^*\{\{\ldots\{\{f_{n_0}+f_{n_0}^*,G_{n_1}\},
G_{n_2}\}\ldots\},G_{n_m}\}
\Big)\\
&=-\Big(
f_{s+1} + f_{s+1}^*+
 \sum_{m\ge1}\frac{1}{m!}\sum_{\substack{n_0+\ldots+n_m=s\\ n_i\ge1}}\{\{\ldots\{\{f_{n_0}+f_{n_0}^*,G_{n_1}\},
G_{n_2}\}\ldots\},G_{n_m}\}\\
&\qquad+\sum_{m\ge2}\frac{1}{m!}\sum_{\substack{n_1+\ldots+n_m=s\\ n_i\ge1}}\{\{\ldots\{\{H_{0},G_{n_1}\},
G_{n_2}\}\ldots\},G_{n_m}\}
\Big)\\
&=-\Big(
f_{s+1} +f_{s+1}^*+
 \sum_{m\ge1}\frac{1}{m!}\sum_{\substack{n_0+\ldots+n_m=s\\ n_i\ge1}}\{\{\ldots\{\{f_{n_0}+f_{n_0}^*,G_{n_1}\},
G_{n_2}\}\ldots\},G_{n_m}\}\\
&\qquad+\sum_{j\ge1}\frac{1}{(j+1)!}\sum_{\substack{n_0+\ldots+n_{j}=s\\ n_i\ge1}}\{\{\ldots\{\{H_{0},G_{n_0}\},
G_{n_1}\}\ldots\},G_{n_{j}}\}
\Big)\,.
\end{aligned}
\end{equation} 
By the inductive hypothesis we have
\begin{equation}\label{induco}
\{H_{0},G_{n_0}\} = -\Big(
f_{n_0}+f_{n_0}^*+
\sum_{l=1}^{n_0-1}c_{l}
\sum_{\substack{h_0+\ldots+h_{l}=n_0\\ h_i\ge1}}\{\{\ldots\{\{f_{h_0}+f_{h_0}^*,G_{h_1}\},
G_{h_2}\}\ldots\},G_{h_l}\}
\Big)\,,
\end{equation}
which inserted into \eqref{Gs+1} gives
\begin{equation}\label{espando}
\begin{aligned}
p\cdot &\del_q G_{s+1} =-\Big(
f_{s+1} + f_{s+1}^* \\
&+\sum_{m\ge1}\frac{1}{m!}\sum_{\substack{n_0+\ldots+n_m=s\\ n_i\ge1}}\{\{\ldots\{\{f_{n_0}+f_{n_0}^*,G_{n_1}\},
G_{n_2}\}\ldots\},G_{n_m}\}\\
&-\sum_{j\ge1}\frac{1}{(j+1)!}\sum_{\substack{n_0+\ldots+n_{j}=s\\ n_i\ge1}}\{\{\ldots\{f_{n_0}+f_{n_0}^*,
G_{n_1}\}\ldots\},G_{n_{j}}\}\\
&-\sum_{j\ge1}\frac{1}{(j+1)!}\sum_{\substack{n_0+\ldots+n_{j}=s\\ n_i\ge1}}
\sum_{l=1}^{n_0-1}c_{l}\\
&\times
\sum_{\substack{h_0+\ldots+h_{l}=n_0\\ h_i\ge1}}
\{\{\ldots\{
\{\{\ldots\{\{f_{h_0}+f_{h_0}^*,G_{h_1}\},
G_{h_2}\}\ldots\},G_{h_l}\}
G_{n_1}\}\ldots\},G_{n_{j}}\}
\Big)\,.
\end{aligned}
\end{equation}

In equation \eqref{espando} above all the indices are strictly positive, and in order to prove \eqref{cm}
we need to collect together the terms with the same number of Poisson brackets, since it is
clear that their coefficient do not depend on the indices $n_i$ or $h_i$.
In the first line the terms with $m$ Poisson brackets appear with a coefficient $1/m!$, in the second line
the terms with $m$ Poisson brackets appear with coefficient $-1/(m+1)!$ and finally in the last line the
terms with $m$ Poisson brackets appear with coefficient
$$
-\sum_{j+l=m}\frac{1}{(j+1)!}c_{l}\,,
$$
therefore the assertion follows.
\EP

\begin{rmk}\label{why}
The cancellation provided by Lemma \ref{lem.cancella} is needed because it allow us
to count the number of summand appearing in \eqref{vera}.
Indeed since all the indices $n_i$ in \eqref{vera} are strictly positive, the number of summand
is equal to the partition function $p(s)$ of the natural number $s$ which,
be the Hardy-Ramanujan asymptotic formula \cite{HR} grows like
\begin{equation}\label{hr}
p(s)\sim\frac{1}{4s \sqrt{3}}{\rm e}^{\pi\sqrt{2s/3}}\ \quad \mbox{ as }\; 
s\to\io\,.
\end{equation}
\end{rmk}


\subsection{Bounds}\label{trinita}

We are now ready to modify the previous construction in order to get the bounds needed for the convergence.

For a function $F=F(q,p)$ and fixed $\s>1$
let us introduce the scale of analytic norms
\begin{equation}\label{lanorma}
\|F\|_{\x}^2=\|F\|_{\x,\s}^2:= \sum_{k\in\ZZZ^2} \avg{k}^{2\s}e^{2\x|k|}\sup_{p\in\calP}|F_{k}(p )|^2,\qquad \x>0
\end{equation}
where we used the standard notation for the Japanese symbol $\avg{k}:=\max\{1,|k|\}$; recall that
since $\s>1$ we have the algebra property
\begin{equation}\label{algebra}
\|FG\|_{\x}\le C_0 \|F\|_{\x} \|G\|_{\x}\,,
\end{equation}
where $C_0$ is some positive constant (depending on $\s$). We omit the index $\s$ in the norm
because it is fixed once and for all.

Set
\begin{equation}\label{numeretti}
\x_s:=\x_0\big(1-\frac{1}{2}\sum_{j=1}^s 2^{-j}\big)\,,
\end{equation}
and note that
\begin{equation}\label{decade}
\x_s> \x_{s+1}\to \frac{\x_0}{2}\,.
\end{equation}

Define
\begin{equation}\label{Fs}
F_{s,k}(p):= \sum_{m=1}^{s-1}c_m |p|^{2m}
\sum_{\substack{n_0+\ldots+n_m=s\\ n_i\ge1}}\{\{\ldots\{\{f_{n_0}+f_{n_0}^*,G^0_{n_1}\},
G^0_{n_2}\}\ldots\},G^0_{n_m}\}\,;
\end{equation}
note that the normalization exponents appearing in \eqref{Fs} are chosen so that
$F_{s,k}(p) = F_{s,k}(p')$ whenever $p,p'$ are parallel, namely $F_{s,k}(p)=\ol{F}_{s,k}(\f_p)$,
where we are using the polar coordinates $(\f_p,r_p)$.
Note also that one has
\begin{equation}\label{realta}
F_{s,-k}(p) = \ol{F}_{s,k}(p)\,.
\end{equation}

\begin{rmk}\label{stocazzo}
The functions $F_{s,k}(p)$ defined in \eqref{Fs} are well defined (and actually analytic) for all $p$ in a neighborhood of $\{H_\e^{(s-1)}=\frac{1}{2}\}\cap\calP$.
\end{rmk}

Let us set
\begin{equation}\label{Gs0}
G^0_{s,k}(p):=\frac{-1}{\ii p\cdot k}\Big(
f_{s,k}+f^*_{s,k}+
F_{s,k}(p)\Big) = \frac{F_{s,k}(p)-F_{s,k}(\tilde{p}_k)}{\ii(p-\hat{p}_{k})\cdot k}\,,
\end{equation} 
}
and note that
\begin{equation}\label{immagino}
G^0_{s,-k}(p) = \ol{G^0}_{s,k}(p)
\end{equation}
so that also the third recursive assuption in \eqref{corsetta} is satisfied.

In particular,
if we set $G^{(s)}=\sum_{m=1}^s \eps^m G^0_m$ 
and $\Phi^{(s)}$ is
the time-$1$ map generated by $G^{(s)}$, by definition we formally get  
\[
H \circ \Phi^{(s)}= \frac{p^2}{2} + \e h_1(p) +\ldots + \e^{s}h_s(p) +\eps^{s+1} R_{s+1}\,,
\]
where
\begin{equation}\label{resto}
\begin{aligned}
\e^{s+1} R_{s+1} &:=\sum_{n=2}^{s}\frac{1}{n!}\{\cdots\{\{ H, 
\sum_{m_1=1}^{s} \e^{m_1} G_{m_1}^0\},\sum_{m_2=1}^{s} \e^{m_2}G_{m_2}^0\}\cdots,
\!\!\!\!\!\!\!\!
\sum_{m=\max\{0,s+1-m_1-\ldots-m_{n-1}\}}^s 
\!\!\!\!
\e^m G_m^0\}\\
&+\sum_{n\ge{s+1}}\frac{1}{n!}\{H,G^{(s)}\}^{(n)}\,,
\end{aligned}
\end{equation}

Note that each summand in \eqref{resto} above has indeed at least a factor $\e^{s+1}$.

By construction we have the following result.

\begin{lemma}\label{stime}
There exist a constant
$C>0$ (depending on $\s$) such that
\begin{subequations}\label{stimoF}
\begin{align}
&\|F_{s}\|_{\x_s}\le \Big( C \| f_1\|_{\x_0}
\Big)^s\,,\label{Effone}\\
&\|f_s\|_{\x_s}, \|\del^h f_s\|_{\x_{s+1}} \le  \Big( C \| f_1\|_{\x_0}
\Big)^s\,,\label{effino}\\
&\|f_s^*\|_{\x_s}, \|\del^h f_s^*\|_{\x_{s+1}} \le  \Big( C \| f_1\|_{\x_0}
\Big)^s\,,\label{effinostar}\\
&\|G^0_s\|_{\x_{s+1}},
\|\del^hG^0_s\|_{\x_{s+1}}\le  \Big( C \| f_1\|_{\x_0}
\Big)^s\,,\label{Gione}
\end{align}
\end{subequations}
where we denoted
$$
\del^h=\del_{p_1}^{\al_1}\del_{p_2}^{\al_2}\del_{q_1}^{\be_1}\del_{q_2}^{\be_2},\qquad 
\al_1+\al_2+\be_1+\be_2=h\,.
$$
\end{lemma}

\prova
We prove the result by induction on $s$;
the case $s=2$ have been considered explicitly in Section \ref{formale}. Assume
inductively the bounds \eqref{stimoF} up to $s$.
The bound \eqref{Effone} follows directly by the definition \eqref{Fs} and Remark \ref{stocazzo}; in particular
the Fourier coefficients will have at most a factor growing
polynomially in $k$ and
we can control such growth with a shrink of the analyticity strip, which provides an
exponential decay $e^{-\x_{s+2} |k|/2}$. Then from  the 
Hardy-Ramanujan formula \eqref{hr} we can bound the number of the summands in \eqref{Fs} and deduce the bound \eqref{Effone}, which 
in turn gives the bound for $f_{s+1}$, $f^*_{s+1}$ and their derivatives. Finally, to get the bound for $G_{s+1}^0$ 
we can reason as done for $G_2^0$. Indeed from the definition \eqref{Gs0} we get
\begin{equation}\label{Gsk}
|G_{s+1,k}^0(p )|=\left|
\frac{F_{s,k}(p)-F_{s,k}(\tilde{p}_k)}{(p-\hat{p}_{k})\cdot k}
\right|
\end{equation}
so that \eqref{Gione} follows
\EP

From Lemma \ref{stime} above we deduce the convergence of the algorithm.

\begin{lemma}\label{stimor}
One has $\|R_{s+1}\|_{\x_{s+1}}\le e^{C\|f_1\|_{\x_0}}$.
\end{lemma}

\prova
Let us factor out the common factor $\e^{s+1}$ in \eqref{resto}. We use the algebra
property \eqref{algebra} and Lemma \ref{stime} above to deduce that
each summand in \eqref{resto} can be bounded by $(C\|f_1\|_{\x_0})^n$,
so that the assertion follows.
\EP

Finally setting
\begin{equation}\nonumber
\begin{aligned}
&G^{(\io)}=\sum_{s=1}^\io \eps^s G^0_s,\\
&f =\sum_{s=1}^\io \eps^s f_s,\\
&f^* =\sum_{s=1}^\io \eps^s f_s^*,
\end{aligned}
\end{equation}
we have the following result.

\begin{lemma}\label{converge}
For all $|\e|<(C \| f_1\|_{\x_0})^{-1}$ one has
\begin{equation}\label{collo}
\| f (\cdot;\e)\|_{\x_0/2}, \| f^* (\cdot;\e)\|_{\x_0/2} ,\; \| G^{(\io)}(\cdot,\cdot;\e)\|_{\x_0/2} <\io.
\end{equation}
\end{lemma}

\prova
One has
$$
\| f (\cdot;\e)\|_{\x_0/2}\le \sum_{s\ge1}\e^s (\| f_s\|_{\x_0/2})^s\le
\sum_{s\ge1}\e^s (\| f_s\|_{\x_s})^s \stackrel{\eqref{effino}}{<}\io.
$$
Same for $f^*$ and $G^{(\io)}$.
\EP

Combining Lemmata \ref{stimor} and \ref{converge} we conclude the proof of Theorem \ref{main}.


\zerarcounters
\section{The $p$-dependent case}\label{azione}

In this section we show that if we allow the correction $f$ to depend also on $p$ then the construction
is much more easy and allows to find local integrability in a whole domain of the phase space.

Now our aim is to explicitely construct a correction 
\begin{equation}\label{pp}
 f(q,p;\e) = \sum_{s\ge1}\e^s f_s(q,p)
 \end{equation}
 and a Hamiltonian $G(q,p;\e)$ such that $H_\e\circ \Phi_G = |p|^2/2$ on an open subset of the phase space,
 where $\Phi_G$ is the time-$1$ flow generated by $G$.
 
We start by choosing $f_1(q,p)$ so that
$$
f_{1,k}(p) = (\ii p\cdot k) g_{1,k}(p),
$$
with $g_{2,k}(p)$  {\it any} function satisfying $|g_{2,k}(p)|\le C Me^{-\x_0|k|}$, so that \eqref{primo} is easily solved.
Clearly such this makes $f_1$ real-on-real.

At order $\e^2$ we still need to solve \eqref{secondo} but now we are on the whole domain
$\PP$
and we allow $f_2$ to depend on $p$. Passing to the Fourier side, we have to solve \eqref{homof};
for all $k$, if we define
\begin{equation}\label{f2p}
f_{2,k}(p):=-\frac{1}{2}\sum_{k_1+k_2=k}
\frac{k_1 f_{1,k_1} k_2 f_{1,k_2}}{(p\cdot k_2)^2} + (\ii p\cdot k)g_{2,k}(p)\,,
\end{equation}
where $g_{2,k}(p)$ is { any} function satisfying $|g_{2,k}(p)|\le C M^2e^{-3\x_0|k|/4}$, then no small divisors appear in \eqref{homof},
and one has simply
$$
G_{2,k}(p) = g_{2,k}(p)\,.
$$
Clearly one can reason exacly in the same way to all orders, by setting
\begin{equation}\label{effesstarp}
f_{s,k}(p):=-\Big(
\sum_{m=1}^{s}\frac{1}{m!}
\sum_{\substack{n_0+\ldots+n_m=s\\ n_0\ge0\\n_1,\ldots,n_m\ge1}}
\{\{\ldots\{\{f_{n_0}+f_{n_0}^*,G^0_{n_1}\},
G^0_{n_2}\}\ldots\},G^0_{n_m}\}\Big) + (ip\cdot k)g_{s,k}(p)\,.
\end{equation}
with $g_{s,k}(p)$ any function satisfying $|g_{2,k}(p)|\le C M^se^{-\x_s|k|}$; in this way the solution of the homological
equation is simply
$$
G_{s,k}(p) = g_{s,k}(p)\,.
$$
Moreover the functions $f_s$ obtained are clearly real-on-real.

This concludes the proof of Theorem \ref{thm-p}. 
\EP


\zerarcounters
\section{Chaotic Behavior}\label{caos}

We now prove that it is possible to choose the potential $f(q,\e)$
in \eqref{ham} so that in the region $\{H_\e=1/2\}\cap\calP^C$ there is at least one orbit that exhibit chaotic behavior.
In order to do so of course we need to be close to a resonance in $p$-space: for instance we may want to be close to
$\ol{p}=(0,1)$.

First of all let $n$ be large enough (to be chosen?) so that
\begin{equation}\label{coseno}
f_{1,(n,0)} = f_{1,(-n,0)} = \frac{1}{2}
\end{equation}
and suppose that
\begin{equation}\label{emmo}
f_1(q) = \cos(nq_1) + \sum_{\substack{k\in\calK \\ |k|<n}}f_{1,k} e^{\ii k\cdot q}\,.
\end{equation}
Of course the high order coefficients of the potential have to be constructed starting from $f_1$ as done in the previous sections.

Let us consider the truncated Hamiltonian
\begin{equation}\label{Htru}
H_1(q,p) = \frac{|p|^2}{2} + \e f_1(q)
\end{equation}
so that the equation of motion are
\begin{equation}\label{bada}
\left\{
\begin{aligned}
\dot{q}_1 &= p_1 \\
\dot{p}_1 &= \e n \sin(nq_1) + \e \sum_{\substack{k\in\calK \\ |k|<n}}\ii k_1f_{1,k} e^{\ii k\cdot q} \\
\dot{q}_2 &= p_2 \\
\dot{p}_2 &= \e \sum_{\substack{k\in\calK \\ |k|<n}}\ii k_2f_{1,k} e^{\ii k\cdot q}\,.
\end{aligned}
\right.
\end{equation}

Let us now introduce the rescaling
\begin{equation}\label{cambio}
p_1 = \sqrt{\e}y_1 \qquad \tau = \sqrt{\e}t
\end{equation}
so that \eqref{bada} can be rewritten as
\begin{equation}\label{bada1}
\left\{
\begin{aligned}
{q}'_1 &= y_1 \\
{y}'_1 &=  n \sin(nq_1) +  \sum_{\substack{k\in\calK \\ |k|<n}}\ii k_1f_{1,k} e^{\ii k\cdot q} \\
{q}'_2 &= \frac{1}{\sqrt{\e}}p_2 \\
{p}'_2 &= \sqrt{\e} \sum_{\substack{k\in\calK \\ |k|<n}}\ii k_2f_{1,k} e^{\ii k\cdot q}\,.
\end{aligned}
\right.
\end{equation}
where $'$ denotes the derivative w.r.t. $\tau$.


\zerarcounters
\section{Final remarks}\label{finale}

The proof of Theorem \ref{main} heavily relies on the dimension; precisely we are able to perform our construction
 because at each step $s$ there is a unique point in $\calP$ such that
\begin{equation}\label{forzaroma}
\begin{aligned}
p\cdot k =0\qquad \qquad \\
\frac{1}{2}|p|^2 + \sum_{k=1}^{s-2}\e^k f_k(q) = \frac{1}{2}\,.
\end{aligned}
\end{equation}

One already sees the problem at the second step when, in order to solve the homological equation \eqref{homof},
one needs to make sure that when $p\cdot k =0$ also the l.h.s. vanishes.

The existence of a unique point in $\calP$ solving \eqref{forzaroma}  allow us to define $f$ as a function
of $q$ only and this is clearly not possible in higher dimension. Indeed if $p\in\RRR^n$ with $n\ge3$ there is a whole
curve $\Gamma$ in $\calP$ solving \eqref{forzaroma}, and thus one seem to be forced to define $f_{2,k}$ as 
a function of $p$, at least for $p\in\Gamma$: of course the same type of argument can be carried at each step.

On the other hand, if we allow $f$ to depend also on $p$, the construction is way easier, as the proof of
Theorem \ref{thm-p} shows, but of course in that case the Hamiltonian \eqref{ham-p} is not associated with a metric.


\end{document}